\documentclass[a4paper,11pt]{amsart} 
\usepackage{amscd}             
\usepackage{amssymb}
\usepackage{amsthm}
\usepackage{epsf} 
\usepackage[T1]{fontenc}            

\newtheorem{thm}{Theorem}[section]   

\newtheorem{lemma}[thm]{Lemma}
\newtheorem{prop}[thm]{Proposition}

\newtheorem{rem}[thm]{Remark}
\renewcommand{\proofname}{Proof}
\def\Cliff{\operatorname{Cliff}}
\def\gon{\operatorname{gon}}
\def\ZZ{{\mathbf Z}}
\def\O{{\mathcal O}}
\def\iso{\simeq}               
\def\+{\oplus}                   
\def\*{\otimes}                  
\def\hpil{\longrightarrow}       

\def\Pic{\operatorname{Pic}}
\def\hs{\hspace{.2cm}}

\begin{document}

\title{Gonality and Clifford Index of curves on $K3$ Surfaces}
\author{Andreas Leopold Knutsen}  
\address{Dept. of Mathematics\\ 
 University of Bergen\\ Johs. Brunsgt 12\\ 5008 Bergen, Norway}
\email{andreask@mi.uib.no}

\begin{abstract}
 We show that every possible value for the Clifford index and gonality
 of a curve of a given genus on a $K3$ surface occurs. 
\end{abstract}

\maketitle

\section{Introduction}
\label{intro}

Curves on $K3$ surfaces have been objects of intense study throughout
the years. In particular, in connection with the famous conjecture of Green
\cite{gre}, Green and Lazarsfeld proved in \cite{gl} that
the Clifford index is the same for all smooth curves in a complete linear system  
on a $K3$ surface. The original conjecture by Harris and Mumford was
that the gonality should be the same for all smooth curves in a linear
system, but this was proved to be false by a counterexample of Donagi
and Morrison \cite{dm}. If the line bundle is ample, then this is
indeed the only counterexample, by a result of Ciliberto and Pareschi
\cite{cp}.

In \cite{elms} curves on $K3$ surfaces are used to construct
exceptional curves of any Clifford dimension. In our recent paper
\cite{kn1} we relate the Clifford index and gonality of curves to the
higher order embeddings of the surface. 

The aim of this note is to show that for $g \geq 3$, genus $g$ curves
of any possible Clifford index and gonality can be found on some $K3$
surface. More precisely, we will show the following statements:

\begin{thm} \label{exthm}
  (a) Let $g$ and $c$ be integers such that $g \geq 4$ and $0 \leq c \leq \lfloor 
  \frac{g-1}{2} \rfloor $. Then there exists a $K3$ surface containing a smooth
  curve of genus $g$ and Clifford index $c$.
 
  (b) Let $g$ and $k$ be integers such that $g \geq 3$ and $2 \leq k \leq \lfloor 
  \frac{g+3}{2} \rfloor $. Then there exists a $K3$ surface containing a smooth
  curve of genus $g$ and gonality $k$.
\end{thm}

Recall that a smooth curve $C$ of genus $g \geq 2$ is said to have {\it
  gonality} $k$ if $C$
posesses a $g^1_k$ but no $g^1_{k-1}$. 

If $A$ is a line bundle
on $C$, then the {\it Clifford index} of $A$ is the integer
\[ \Cliff A = \deg A - 2(h^0 (A) -1). \]
When $g \geq 4$, the {\it Clifford index of $C$} itself is defined as 
\[ \Cliff C = \min \{ \Cliff A \hspace{.05in} | \hspace{.05in} h^0 (A) \geq 2, h^1 (A) \geq 2 \}. \]

Note that we always have $\Cliff C \leq \gon C-2$. Furthermore, we get from Brill-Noether theory that the gonality of $C$ satisfies $\gon C \leq \lfloor 
  \frac{g+3}{2} \rfloor $, whence $\Cliff C \leq \lfloor 
  \frac{g-1}{2} \rfloor $. For the general curve of genus
  $g$, we have $\gon C = \lfloor 
  \frac{g+3}{2} \rfloor $ and $\Cliff C = \lfloor \frac{g-1}{2}
  \rfloor $. 

We say that a line bundle $A$ on $C$ {\it contributes to the Clifford
  index of $C$} if $h^0 (A)\geq 2$ and $h^1 (A) \geq 2$ and that it {\it computes 
the Clifford index of $C$} if in addition $\Cliff C = \Cliff A$. 

It is classically known that for a fixed genus, curves of any possible
gonality occur. By Ballico's result \cite{bal}, the same is true for
the Clifford index. Our result then gives a new proof of these facts,
with the additional information that one can always construct such
curves on $K3$ surfaces.

The central tool in this paper is the already mentioned famous result of Green and Lazarsfeld \cite{gl}, which states
that if $L$ is a base point free line bundle on a $K3$ surface $S$, then
$\Cliff C$ is constant for all smooth irreducible $C \in |L|$, and if
$\Cliff C < \lfloor \frac{g-1}{2} \rfloor $, then there exists a line
bundle $M$ on $S$ such that $M_C := M \* \O _C$ computes the Clifford index of
$C$ for all smooth irreducible $C \in |L|$. (Note that since $(L-M) \*
\O _C \iso \omega _C \* {M_C}^{-1}$, the result is symmetric in $M$ and
$L-M$.)

Other tools are classical results on linear systems on $K3$ surfaces
\cite{S-D} and lattice theory of $K3$ surfaces \cite{morrison}.

For other results on the Clifford index and gonality of smooth curves
on $K3$ surfaces we refer to \cite{cp}, \cite{elms} and \cite{kn1}.

A {\it curve} is always reduced and irreducible.

The author would like to thank Gavril Farkas for interesting conversations.

\begin{rem}
  {\rm We adopt the same convention as in \cite{gl}, and set $\Cliff C =0$ for
  $C$ of genus $2$ or hyperelliptic of genus $3$, and $\Cliff C =1$ for $C$
  non-hyperelliptic of genus $3$. The result of Green and Lazarsfeld is then
  still valid.} 
\end{rem}

\section{Proof of the main theorem}
\label{proof}

The theorem is an immediate consequence of the following:

\begin{prop} \label{exprop}
  Let $d$ and $g$ be integers such that $g \geq 3$ and $2 \leq d \leq \lfloor \frac{g+3}{2} \rfloor $, and let $S$ be a $K3$
  surface with $\Pic S = \ZZ L \+ \ZZ E$, where $L^2 = 2(g-1)$, $E.L = d$ and $E^2 = 0$.
  Then $L$ is base point free and for any smooth curve $C \in |L|$ we
  have 
  \[ \Cliff C = d-2 \leq \lfloor \frac{g-1}{2} \rfloor \]
    and $\Cliff C$ is computed by the pencil $\O_C(E)$.
\end{prop}

To prove this proposition, we first need the following basic existence result:

\begin{lemma} \label{help}
  Let $g \geq 3$ and $d \geq 2$ be integers. Then there exists a $K3$
  surface $S$ with $\Pic S = \ZZ L \+ \ZZ E$, such that $L$ is base point free and $E$ is a smooth curve, $L^2 = 2(g-1)$, 
  $E.L = d$ and $E^2 = 0$.
\end{lemma}

\begin{proof}
  By \cite[Prop. 4.2]{knut}, we can find a $K3$
  surface $S$ with $\Pic S = \ZZ L \+ \ZZ E$, with intersection matrix 
\[  \left[ 
  \begin{array}{cc}
  L^2  &  L.E   \\ 
  E.L  &  E^2 
    \end{array} \right]  = 
    \left[
  \begin{array}{cc}
  2(g-1)  &  d   \\ 
  d   &  0      
    \end{array} \right]     \] 
  and such that $L$ is nef. This is a consequence of the lattice
  theory in \cite{morrison}, which again follows from the surjectivity of the
  period map. If $L$ is not base point free, there
  exists by standard results on linear systems on $K3$ surfaces (see
  e.g. \cite[2.7]{S-D} or \cite[Thm. 1.1]{kn1}) a curve 
  $B$ such that $B^2=0$ and $B.L =1$. An easy calculation, writing $B
  \sim xL + yE$ for two integers $x$ and $y$, 
  shows that this is 
  impossible. By \cite[Proposition 4.4]{knut}, we have that  $|E|$
  contains a smooth curve. 
\end{proof}

\renewcommand{\proofname}{Proof of Proposition {\rm \ref{exprop}}}
  
  \begin{proof} 
  Let $S$, $L$ and $E$ be as in Lemma \ref{help}, with 
  $d \leq \lfloor \frac{g+3}{2} \rfloor $. Note that 
  since $E$ is irreducible, we have $h^1(E) = 0$. Furthermore,
  Riemann-Roch gives  $\chi (L-E) = \frac{1}{2} (L-E)^2+2 =g+1-d \geq
  \frac{g-1}{2} \geq 1$, whence by Serre duality either $h^0(L-E) \geq
  1$ or $h^0(E-L) \geq 1$ (but not both simultaneously). Since
  $(E-L).L =d-2g+2 <0$, we have $h^0(L-E) \geq
  1$. If $g \geq 4$, we even have $h^0(L-E) \geq 2$. 

  By the
  cohomology of the short exact sequence
\[ 0 \hpil \O_S (E-L) \hpil \O_S (E) \hpil \O_{C  } (E) \hpil 0, \]
  where $C$ is any smooth curve in $|L|$, we find that $h^0 (\O_{C  } (E)) \geq h^0(E) =2$ and 
  $h^1 (\O_{C  } (E))= h^0(L-E)$, so $\O_{C  } (E)$ contributes to the 
Clifford index of
$C$ (if $g \geq 4$) and 
\[ c:= \Cliff C \leq \Cliff \O_{C  } (E) \leq E.L - E^2 - 2 = d-2 < \lfloor
  \frac{g-1}{2} \rfloor. \]

If $c = \lfloor
  \frac{g-1}{2} \rfloor$, there is nothing more to prove. 

If $c < \lfloor \frac{g-1}{2} \rfloor $, then by
\cite[(2.3)]{mar} and \cite[Lemma 8.3]{kn1}
there exists an effective divisor $D$ on $S$ satisfying 
\begin{eqnarray}
 \label{eq:1} h^0 (\O_C(D))=h^0(D) \geq 2, \hs h^1(\O_C(D))=h^0(L-D) \geq 2, \\
 \label{eq:2} h^1(D)=h^1(L-D)=0  \hs \mbox{ and } \\
 \label{eq:3} c = \Cliff \O_{C} (D) = \Cliff \O_{C} (L-D)= D.L - D^2 - 2
\end{eqnarray}
(this is a consequence of the result of Green and Lazarsfeld \cite{gl} mentioned above).

Since both $D >0$ and $L-D >0$ and $E$ is nef, we must have
\[ D.E \geq 0 \mbox{     and     } (L-D).E \geq 0. \]
Writing $D \sim xL + yE$ this is equivalent to
\[ dx  \geq 0 \mbox{     and     }  d(1-x)  \geq 0, \]
which gives $x=0$ or $1$. These two cases give, respectively, $D=yE$ and
$L-D =yE$. Since $h^1(yE)=y-1$ by \cite[Prop. 2.6]{S-D}, it follows by
(\ref{eq:2}) that $y=1$ and $D \sim E$ or $L-D \sim E$. Hence $c = \Cliff
\O_C(E)= E.L - E^2 - 2 = d-2$ by (\ref{eq:3}). By Riemann-Roch and (\ref{eq:1}) we
also have $h^0 (\O_C(E))=h^0(E)= \frac{1}{2}E^2+2=2$, so $\O_C(E)$ is
a pencil. 
\end{proof}

\renewcommand{\proofname}{Proof}

This concludes the proof of Theorem \ref{exthm}.

\vspace{.5cm}

\end{document}